\documentclass{article}
\usepackage{amssymb}
\usepackage{amsmath}
\usepackage{amsthm}
\usepackage{mathrsfs}
\usepackage{amsfonts}
\begin{document}
\topmargin= -.2in \baselineskip=15pt
\newtheorem{theorem}{Theorem}[section]
\newtheorem{proposition}[theorem]{Proposition}
\newtheorem{lemma}[theorem]{Lemma}
\newtheorem{corollary}[theorem]{Corollary}
\newtheorem{conjecture}[theorem]{Conjecture}
\theoremstyle{remark}
\newtheorem{remark}[theorem]{Remark}

\title {$\ell$-adic Realization of Some Aspects of Landau-Ginzburg $B$-models\thanks{I would
like to thank the referee for many suggestions. This research is
supported by NSFC.}}

\author {Lei Fu\\
{\small Yau Mathematical Sciences Center, Tsinghua University,
Beijing, China}\\
{\small leifu@mail.tsinghua.edu.cn}}
\date{}
\maketitle

The Landau-Ginzburg $B$-model for a germ of a holomorphic function
with an isolated critical point is constructed by K. Saito \cite{KS}
and finished by M. Saito \cite{MS}. Douai and Sabbah construct the
Landau-Ginzburg $B$-models for some Laurent polynomials \cite{D1,
D2, DS}. The construction relies on analytic procedures, and one can
not expect it can be done by purely algebraic method. In this note,
we  work out the $\ell$-adic realization of the algebraic part of
the construction. In \S 1, we define Frobenius type structures. One
can consult \cite{HM,S1,S2} for details. In \S 2, we sketch the
construction of the Landau-Ginzburg $B$-models for Laurent
polynomials. Details can be found in \cite{D1, D2, DS, HM}. In \S 3,
we study the $\ell$-adic counterpart of the construction in \S 2.

\section{Frobenius type structures}

Let $D$ be the germ of $\mathbb C$ at $0$,  let $X$ be a germ of complex manifold, and
let $\mathcal F$ be a trivial holomorphic vector bundle on $D\times X$. Denote the
$\mathcal O_{D\times X}$-module of holomorphic sections of $\mathcal F$ also by $\mathcal F$.
Let  $$\nabla:\mathcal F|_{(D-\{0\})\times X}\to \mathcal
F|_{(D-\{0\})\times X}\otimes \Omega^1_{(D-\{0\})\times X}$$ be an integrable connection. We say $\nabla$ has a
pole of \emph{Poincar\'e rank} $\leq m$ along the divisor $X\times 0$  if
$$\nabla(\mathcal F) \subset \mathcal F\otimes\frac{1}{t^m} \Big(\sum_i\mathcal O_{D\times X}dx_i+
\mathcal O_{D\times X}\frac{dt}{t}\Big),$$ where $t$ is the coordinate for $D$, and $(x_i)$ the coordinate for $X$.

\medskip
\noindent{\bf Poincar\'e rank 0 case:}

If the Poincar\'e rank is $0$, we say $(\mathcal F,\nabla)$ has a
\emph{logarithmic pole} along $0\times X$. Fix a global basis for
$\mathcal F$, and write the connection matrix as
$$A=\sum_i
\Omega_i(t,x)dx_i+\Omega(t,x)\frac{dt}{t}$$ for some matrices of
holomorphic functions $\Omega_i(t,x)$ and $\Omega(t,x)$.
We have
\begin{eqnarray*}
dA+A\wedge A&=&\sum_i \Big(\frac{1}{t}\big(-\frac{\partial \Omega}{\partial x_i}+[\Omega,\Omega_i]\big) + \frac{\partial \Omega_i}{\partial t}\Big)
dtdx_i + \sum_{i<j} \Big(\frac{\partial \Omega_j}{\partial x_i}-\frac{\partial \Omega_i}{\partial x_j} +[\Omega_i,\Omega_j] \Big)dx_idx_j.
\end{eqnarray*}
Since
$dA+A\wedge A=0,$ we have
\begin{eqnarray*}
\frac{1}{t}\Big(-\frac{\partial \Omega}{\partial x_i}+[\Omega,\Omega_i]\Big) + \frac{\partial \Omega_i}{\partial t}&=&0,\\
\frac{\partial \Omega_j}{\partial x_i}-\frac{\partial \Omega_i}{\partial x_j} +[\Omega_i,\Omega_j]&=&0.
\end{eqnarray*}
It follows that
\begin{eqnarray*}
&&\Big(\frac{\partial \Omega}{\partial x_i}-
[\Omega,\Omega_i]\Big)|_{0\times X}=0,\\
&& d(\sum_i \Omega_i(0,x)dx_i)+(\sum_i \Omega_i(0,x)dx_i)\wedge
(\sum_i \Omega_i(0,x)dx_i)=0.
\end{eqnarray*}
The second equation shows that $\sum_i \Omega_i(0,x) dx_i$ defines
an integrable connection $\triangledown=\nabla|_{0\times X}$ on
$\mathcal F|_{0\times X}$. The first equation shows that
$\Omega(0,x)$ defines a horizontal endomorphism
$R_0=\mathrm{Res}_0(\nabla)$ of $(\mathcal F|_{0\times X},
\triangledown)$, which we call \emph{the residue of} $\nabla$. We
summarize the above data as
\begin{eqnarray*}
\triangledown=\nabla|_{0\times X},\quad
R_0=\mathrm{Res}_0(\nabla)=\nabla_{t\partial_t}|_{0\times X},\quad
\triangledown\triangledown =0,\quad \triangledown(R_0)=0.
\end{eqnarray*}

\medskip
\noindent{\bf Poincar\'e rank $\leq 1$ case:}

The connection matrix is of the form
$$A=\frac{1}{t}\Big(\sum_i\Omega_i(t,x)dx_i+\Omega(t,x)\frac{dt}{t}\Big)$$
for some matrices of holomorphic functions $\Omega_i(t,x)$ and
$\Omega(t,x)$. We have
\begin{eqnarray*}
&&dA+A\wedge A\\
&=& \sum_i\Big(\frac{1}{t^3}
[\Omega,\Omega_i]-\frac{1}{t^2}\big(\frac{\partial \Omega}{\partial x_i}+\Omega_i\big) +\frac{1}{t}\frac{\partial \Omega_i}{\partial t}\Big)dtdx_i \\
&&\quad \sum_{i<j}\Big(\frac{1}{t^2}[\Omega_i,\Omega_j]+\frac{1}{t}\big(\frac{\partial \Omega_j}{\partial x_i}-\frac{\partial \Omega_i}{\partial x_j}   \big)\Big) dx_idx_j
\end{eqnarray*}
From the equation $dA+A\wedge A=0$, one
deduces that
$$[\Omega, \Omega_i]|_{0\times X}=0,\quad [\Omega_i, \Omega_j]|_{0\times
X}=0.$$ Let
$$\Phi=(t\nabla)|_{0\times X}=\sum\Omega_i(0,x)dx_i:\mathcal F|_{0\times
X}\to \mathcal F|_{0\times X}\otimes \Omega^1_X.$$ The equation
$[\Omega_i, \Omega_j]|_{0\times X}=0$ shows that
$$\Phi\wedge \Phi=0,$$ that is, $\Phi$ is a Higgs field on $\mathcal
F|_{0\times X}$.  Let
$$R_0=(t^2\nabla_{\partial_t})|_{0\times
X}=\Omega(0,x)\in \mathrm{End}(\mathcal F|_{0\times X}).$$ The
equation $[\Omega, \Omega_i]|_{0\times X}=0$ shows that
$$[R_0,\Phi]=0.$$

\medskip
Denote the standard coordinate on $\mathbb A^1=\mathbb
P^1-\{\infty\}$ by $t$ and let $s=\frac{1}{t}$. Let $\mathcal E$ be
a trivial vector bundle on $X$. Denote the sheaf of holomorphic
sections of $\mathcal E$ by the same notation. Let $\pi:\mathbb
P^1\times X\to X$ be the projection. Suppose we have an integrable
connection $\nabla$ on $\pi^\ast \mathcal E$ with a logarithmic pole
along $\infty\times X$, a pole of Poincar\'e rank $\leq 1$ along
$0\times X$, and holomorphic elsewhere. Then $\mathcal E\cong
(\pi^\ast \mathcal E)|_{\infty\times X}$ is endowed with an
integrable connection
$$\triangledown =\nabla|_{\infty\times X}$$
and a horizontal endomorphism
$$R_\infty=\mathrm{Res}_\infty(\nabla),$$
and $\mathcal E\cong (\pi^\ast\mathcal E)|_{0\times X}$
is also endowed with a Higgs field
$$\Phi=(t\nabla)|_{0\times X}$$
and an endomorphism $$R_0=(t^2\nabla_{\partial_t})|_{0\times X}$$
commuting with $\Phi$. The connection $(\nabla, \pi^\ast\mathcal E)$
can also be constructed from the tuple $(\mathcal E, \triangledown,
R_0,R_\infty,\Phi)$. This gives rise the to the so-called Frobenius
type structure. More precisely, a \emph{Frobenius type structure}
(without a metric) on $X$ is a tuple $(\mathcal E, \triangledown,
R_0, R_\infty, \Phi)$ such that $\mathcal E$ is free $\mathcal
O_X$-module of finite rank, $R_0, R_\infty\in \mathrm{End}_{\mathcal
O_X}(\mathcal E)$, $\Phi:\mathcal E\to \mathcal E\otimes_{\mathcal
O_X}\Omega_X^1$ is a Higgs field, $\triangledown: \mathcal E\to
\mathcal E\otimes_{\mathcal O_X}\Omega_X^1$ is an integrable
connection, and we require the following condition holds. Let
$\pi:\mathbb P^1\times X\to X$ be the projection. We require that
the connection
\begin{eqnarray}
\begin{array}{rcl}
\nabla&=&\pi^\ast \triangledown +\frac{\pi^\ast
\Phi}{t}+\Big(\frac{R_0}{t}-R_\infty\Big)\frac{dt}{t}\\
&=&\pi^\ast \triangledown +s {\pi^\ast \Phi}+(-s
{R_0}+R_\infty)\frac{ds}{s} \end{array}
\end{eqnarray} on $\pi^\ast
\mathcal E$ is integrable. Note that $\nabla$ has a logarithmic pole
along $\infty\times X$, a pole of Poincar\'e rank $\leq 1$ along
$0\times X$, and holomorphic elsewhere. The condition
$\nabla\nabla=0$ is equivalent to the conditions
\begin{eqnarray*}
&&\triangledown \triangledown =0, \quad \triangledown(R_\infty)=0,
\quad \Phi\wedge \Phi=0, \quad[\Phi, R_0]=0, \\
&& \triangledown(\Phi)=0, \quad
\triangledown(R_0)+\Phi=[\Phi,R_\infty].
\end{eqnarray*}
Indeed, fix a global basis for $\mathcal E$, and let $A$, $\Phi$,
$R_0$, $R_\infty$ be the matrix for the connection $\triangledown$,
the Higgs field $\Phi$, the endomorphisms $R_0$ and $R_\infty$,
respectively. Note that $A$ and $\Phi$ are matrix of holomorphic
$1$-forms on $X$, and $R_0$ and $R_\infty$ are matrix of holomorphic
functions on $X$. The connection matrix for $\nabla$ is
$A'=A+\frac{\Phi}{t}+\Big(\frac{R_0}{t}-R_\infty\Big)\frac{d t}{t}$.
The expression for $dA'+A'\wedge A'$ is
 \begin{eqnarray*}
&& (dA+A\wedge A)
+\Big(\frac{1}{t^3}[\Phi,R_0]+\frac{1}{t^2}\big(dR_0+[A,R_0]+\Phi-[\Phi,R_\infty]
 \big)-\frac{1}{t}\big(dR_\infty +[A,R_\infty] \big)  \Big) dt \\
 &&\quad +\frac{1}{t^2}\Phi\wedge \Phi+ \frac{1}{t}\Big(d\Phi+A\Phi+\Phi A \Big).
 \end{eqnarray*}
Setting it equals to $0$, we get the relations above.

Any meromorphic integrable connection $\nabla$ on a trivial vector
bundle $\pi^\ast\mathcal E$ over $\mathbb P^1\times X$ with a
logarithmic pole along $\infty\times X$, a pole of Poincar\'e $\leq
1$ along $0\times X$, and holomorphic elsewhere is called a trTLE
structure by Hertling. One can show any trTLE structure is of the
form (1) and hence gives rise to a Frobenius type structure without
a metric.  Confer \cite[Theorem 4.2]{HM}. By abuse of notation, we
also call a trTLE structure a Frobenius type structure.

\medskip
\noindent {\bf Birkhoff problem}: Let $D$ be a disc, and let $(\mathcal F, \nabla)$ be a
trivial holomorphic bundle on $D$ equipped with an
integrable connection with a pole of Poincar\'e rank $\leq 1$ at $0$.
Find a pair $(\widetilde{\mathcal F}, \widetilde \nabla)$ such that $\widetilde {\mathcal F}$ is a trivial bundle on
$\mathbb P^1$,  $\widetilde \nabla$ is an integrable meromorphic connection
with logarithmic pole at $\infty$ and holomorphic outside $\{0,\infty\}$, and
$(\widetilde{\mathcal F}, \widetilde \nabla)|_D\cong ({\mathcal F}, \nabla)$.

\begin{proposition}[Birkhoff problem for a family] Let $D$ be a disc, $(X,x_0)$
a germ of complex manifold, and $(\mathcal F,\nabla)$ a trivial
holomorphic bundle on $D\times X$ equipped with an integrable
meromorphic connection of Poincar\'e rank $\leq 1$ along $0\times X$.
Suppose we can solve the Birkhoff problem for $(\mathcal
F,\nabla)|_{D\times \{x_0\}}$. Then there exists a unique pair $(\widetilde {\mathcal F}, \widetilde \nabla)$
of a trivial holomorphic vector bundle $\widetilde {\mathcal F}$ on
$\mathbb P^1\times X$ equipped with an integrable meromorphic
connection $\widetilde \nabla$ with logarithmic pole along $\infty\times X$ and holomorphic outside $\{0,\infty\}\times X$,
such that $(\widetilde {\mathcal F}, \widetilde \nabla)|_{\mathbb P^1\times \{x_0\}} $ is the given solution of the Birkhoff
problem, and $(\widetilde {\mathcal F}, \widetilde \nabla)|_{D \times X}\cong
(\mathcal F,\nabla)$. We
thus get a Frobenius type structure.  We have $(\widetilde{\mathcal F},\widetilde\nabla)|_{D_\infty\times X}\cong p^\ast
((\widetilde {\mathcal F},\widetilde \nabla)|_{D_\infty\times\{x_0\}})$, where $D_\infty$ is a disc centered at $\infty$, and
$p:D_\infty\times X\to D_\infty\times\{x_0\}$ is the projection.
\end{proposition}

\noindent{\bf Algebraic version of the Birkhoff problem}. Let $(G,\nabla)$ be a
$\mathbb C[t, t^{-1}]$-module equipped with a connection
having poles only at $0,\infty$ with regular singularity at
$\infty$,  and let $G_0$ be a free $\mathbb C[t]$-submodule of $G$
with Poincar\'e rank $\leq 1$ such that $G_0\otimes_{\mathbb
C[t]}\mathbb C[t,t^{-1}] =G$. Find a free
$\mathbb C[t^{-1}]$-submodule $G_\infty$ of $G$ which is
logarithmic such that
$$G_0\otimes_{\mathbb C[t]}\mathbb C[t,t^{-1}]=G_\infty\otimes_{\mathbb C[t^{-1}]}\mathbb C[t,t^{-1}]
=G$$ and such that the vector bundle on $\mathbb P^1$ obtained by gluing $G_0$ and $G_\infty$ is free.

\medskip
Suppose the monodromy of $(G, \nabla)$ at $\infty$ is
quasi-unipotent, and let $s=\frac{1}{t}$. Let $V_\bullet\mathbb
C[s]\langle{\partial_s}\rangle$ be the increasing filtration of the
Weyl algebra $C[s]\langle{\partial_s}\rangle$ defined by
\begin{eqnarray*}
V_{-k}\mathbb C[s]\langle{\partial_s}\rangle&=& s^k\mathbb
C[s]\langle s {\partial_s}\rangle\quad \hbox {for }k\geq
0,\\
V_{k}\mathbb C[s]\langle{\partial_s}\rangle&=& V_{k-1} \mathbb
C[s]\langle{\partial_s}\rangle+{\partial_s}  V_{k-1} \mathbb
C[s]\langle{\partial_s}\rangle \quad \hbox {for }k\geq 1.
\end{eqnarray*}
There exists a unique increasing exhaustive filtration $V_\bullet G$
of $G$, indexed by a union of a finite number of subsets
$\alpha+\mathbb Z$ $(\alpha\in\mathbb Q)$ satisfying the following
condition:

(a) For every $\alpha$, the filtration $V_{\alpha+\mathbb Z}G$ is
good relative to $V_\bullet\mathbb C[s]\langle{\partial_s}\rangle$ ;

(b) For every $\beta\in \mathbb Q$, $s\partial_s+\beta$ is nilpotent
on $\mathrm{Gr}_\beta^V(G)=V_\beta G/V_{<\beta} G$. Denote this
nilpotent endomorphism on $\mathrm{Gr}_\beta^V(G)$ by $N$.

One can show each $V_\beta G$ is a free $\mathbb C[s]$-module,
$\mathbb C[s,s^{-1}]\otimes_{\mathbb C[s]} V_\beta G\cong G$, and
the connection has logarithmic pole on $V_\beta G$ at $s=0$.
Consider also the increasing filtration $G_\bullet$ of $G$ defined
by $G_k= t^{-k}G_0$. This filtration induces a filtration
$G_\bullet(\mathrm{Gr}_\beta^V(G))$ on $\mathrm{Gr}_\beta^V(G)$ for
each $\beta$. Let
$$H_\infty=\bigoplus_{\beta\in [0,1)}\mathrm{Gr}_\beta^V (G).$$ It
is the nearby cycle of $(G,\nabla)$ at $\infty$.

\begin{theorem}[M. Saito's criterion] Suppose there exists
a mixed Hodge structure on the nearby cycle $H_\infty$ so that the
Hodge filtration is $G_\bullet H_\infty$, and the weight filtration
is the monodromy filtration $M_\bullet H_\infty$ of $N$. Then we can
solve the Birkhoff problem.
\end{theorem}

We propose the following problem:

\medskip
\noindent {\bf $\ell$-adic version of the Birkhoff problem.} Let $\mathbb F_q$
be a finite field with $q$ elements of characteristic $p$, let $\ell$ be a prime number distinct from $p$, let $\eta_0$
be the generic point of the henselization of
$\mathbb P^1_{\mathbb F_q}$ at $0$, and let
$\rho: \mathrm{Gal}(\overline \eta_0/\eta_0)\to
\mathrm{GL}(n,\overline{\mathbb Q}_\ell)$ be a $\overline{\mathbb
Q}_\ell$-representation. Find conditions on $\rho$ under
which there exists a lisse punctually pure $\overline{\mathbb
Q}_\ell$-sheaf $\mathcal F$ on $\mathbb P^1-\{0,\infty\}$ tamely
ramified at $\infty$ so that $\mathcal F|_{\eta_0}$ corresponds to
the given representation $\rho$.

\medskip
Let $\eta_\infty$ be the generic point of the henselization of
$\mathbb P^1_{\mathbb F_q}$ at $\infty$. Then $\mathcal
F_{\bar\eta_\infty}$ considered as a representation of
$\mathrm{Gal}(\bar\eta_\infty/\eta_\infty)$ is the nearby cycle of
$\mathcal F$ at $\infty$. Since $\mathcal F$ is pure, the monodromy
filtration on $\mathcal F_{\bar\eta_\infty}$ is the weight
filtration (up to a shift) by \cite[1.8.4]{D}, and this corresponds
exactly to the condition of Saito's criterion. The condition that
$\tilde{\mathcal F}$ is a trivial bundle over $\mathbb P^1$ in the
classical Birkhoff problem is replaced by the condition that
$\mathcal F$ is a lisse punctually pure sheaf on $\mathbb
P^1-\{0,\infty\}$ in the $\ell$-adic version.

\medskip
Finally a \emph{Frobenius type structure with a metric} on $X$ is a
tuple $(\mathcal E, \triangledown, R_0, R_\infty, \Phi,g)$ such that
$(\mathcal E, \triangledown, R_0, R_\infty, \Phi)$ is a Frobenius
type structure defined above, and $g: \mathcal E\times\mathcal E\to
\mathcal O_X$ is a symmetric non-degenerate $\triangledown$-flat
pairing such that the pairing
$$G:\pi^\ast \mathcal E \times a^\ast \pi^\ast \mathcal E\to \mathcal
O_{\mathbb P^1\times X}$$ induced by $g$ is $\nabla$-flat, where
$a:\mathbb P^1\to\mathbb P^1$ is the morphism $a(t)=-t$. This is
equivalent to saying that
$$g(\Phi_v(a), b)=g(a,\Phi_v(b)),\quad g(R_0(a), b)=g(a,R_0(b)),\quad
g(R_\infty(a), b)=-g(a,R_\infty(b))$$ for any tangent vector $v$ of
$X$ and any sections $a$ and $b$ of $\mathcal E$. Frobenius type
structures with a metric correspond to trTLEP structures of Hertling
(\cite[Theorem 4.2]{HM}).

\section{Frobenius type structures associated to a subdiagram deformation}

Let $\mathbf{w}_j=(w_{1j},\ldots, w_{nj})\in\mathbb Z^n$
$(j=1,\ldots, N)$, let $f=\sum_{j=1}^N a_j t_1^{w_{1j}}\cdots
t_n^{w_{nj}}$  be a Laurent polynomial with nonzero coefficients
$a_j\in \mathbb C$,  and let $\Delta$ be the \emph{Newton
polyhedron} of $f$ at $\infty$, that is, the convex hull of the set
$\{0, \mathbf{w}_1,\ldots, w_N\}$ in $\mathbb R^n$. We say $f$ is
convenient if $0$ lies in the interior of $\Delta$. We say $f$ is
\emph{non-degenerate} if for any face $\sigma$ of $\Delta$ not
containing $0$, the equations
$$\frac{\partial f_\sigma}{\partial t_1}=\cdots=\frac{\partial f_\sigma}{\partial
t_n}=0$$ define an empty subscheme in $\mathbb
G_m^n=\mathrm{Spec}\,\mathbb C[t_1^{\pm 1},\ldots, t_n^{\pm 1}]$,
where $f_\sigma=\sum_{\mathbf w_j\in \sigma} a_jt_1^{w_{1j}}\cdots
t_n^{w_{nj}}$.

Let $g_1(t_1,\ldots, t_n), \ldots, g_m(t_1,\ldots, t_n)$ be a
family of Laurent polynomials. Consider the deformation
\begin{eqnarray*}
F_x(t)=f(t_1,\ldots, t_n)+x_1 g_1(t_1,\ldots, t_n)+\cdots + x_m
g_m(t_1,\ldots, t_n)
\end{eqnarray*} of $f$.  We say $F$ is a
\emph{subdiagram deformation} of $f$ if all exponents of monomials
with nonzero coefficients in $g_1,\ldots, g_m$ lie in the interior
of $\Delta$.  We say $F$ is the
\emph{universal unfolding} of $f$ if the images of $g_1,\ldots, g_m$
in the Jacobian quotient ring $\mathbb C[t_1^{\pm 1},\ldots,
t_n^{\pm 1}]/\Big(\frac{\partial f}{\partial t_1},\ldots,
\frac{\partial f}{\partial t_n}\Big)$ form a basis.

Suppose $f$ is convenient and non-degenerate, and $F$ is a subdiagram deformation.
Consider the twisted algebraic de Rham complex
$$(\Omega^\bullet_{{\mathbb G}_m^n\times\mathbb A^m/\mathbb A^m}[\tau,\tau^{-1}],
e^{\tau F}\circ d\circ e^{-\tau F}).$$ We have
$$ e^{\tau F}\circ d\circ e^{-\tau F}(\omega)=d
\omega-\tau dF\wedge \omega.$$  Let
\begin{eqnarray*}
G&=&\frac{\Omega^n_{{\mathbb G}_m^n\times\mathbb A^m/\mathbb A^m}[\tau,
\tau^{-1}]}{\Big(d-\tau dF\wedge\Big)\Omega^{n-1}_{{\mathbb G}_m^n\times\mathbb A^m/\mathbb A^m}[\tau, \tau^{-1}]},\\
G^{(o)} &=&\Omega^n_{{\mathbb G}_m^n}[\tau,\tau^{-1}]\large/\Big(d
-\tau df\wedge\Big)\Omega^{n-1}_{\mathbb T^n}[\tau,\tau^{-1}].
\end{eqnarray*} One can show $G$ is a free $\mathbb
C[x_1,\ldots, x_m][\tau, \tau^{-1}]$-module of rank $n!\mathrm{vol}(\Delta)$,
and hence defines a
trivial vector bundle over $(\mathbb P^1-\{0,\infty\})\times \mathbb A^m$.
Define a connection $\nabla$ on $G$ and $G^{(o)}$ by
\begin{eqnarray*}
\nabla_{\partial_{x_j}}=e^{\tau F}\circ{\partial_{x_j}}\circ
e^{-\tau F},\quad \nabla_{\partial_\tau}=e^{\tau F}\circ {\partial
_\tau}\circ e^{-\tau F}.
\end{eqnarray*}
We have
\begin{eqnarray*}
\nabla_{\partial_{x_j}}(\omega)=\frac{\partial\omega}{\partial
x_j}-\tau \frac{\partial F}{\partial x_j} \omega,\quad
\nabla_{\partial_\tau}(\omega)=\frac{\partial\omega}{\partial
\tau}-F \omega.
\end{eqnarray*}
One can see that $(G,\nabla)$ has regular singularity at $\tau=0$.
Let $\theta=\frac{1}{\tau}$. Set
\begin{eqnarray*}
G_0 &=&\frac{\Omega^n_{{\mathbb G}_m^n\times\mathbb A^m/\mathbb
A^m}[\theta]}{\Big(\theta
d- dF\wedge\Big)\Omega^{n-1}_{{\mathbb G}_m^n\times\mathbb A^m/\mathbb A^m}[\theta]},\\
G_0^{(o)} &=&\Omega^n_{{\mathbb G}_m^n}[\theta]\large/\Big(\theta
d-df\wedge\Big)\Omega^{n-1}_{{\mathbb G}_m^n}[\theta].
\end{eqnarray*}
$G_0$ is a free $\mathbb C[x_1,\ldots, x_m][\theta]$-module. It
defines a trivial vector bundle over $(\mathbb
P^1-\{\theta=\infty\})\times \mathbb A^m$. It is a lattice of $G$,
and $\nabla$ defines a meromorphic connection on $G_0$ with
Poincar\'e rank $\leq 1$ along the divisor $\theta=0$. We call $G_0$
the \emph{Brieskorn lattice} associated to the subdiagram
deformation $F$. Using Saito's criterion, Douai and Sabbah \cite{DS}
prove that the Birkhoff problem is solvable for the pairs $(G,G_0)$
(family version) and $(G^{(o)}, G_0^{(o)})$.

\medskip
Next suppose $F$ is a universal unfolding. The definition of the pair $(G,G_0)$
requires some analytic procedure due to the disappearance at infinity of critical points of $F_x(t)$ as
$x\to 0$.
Roughly speaking, $G$ is the Fourier transform
of the Gauss-Manin system for the family $F_x(t)$. There exists a neighborhood $X$ of $0$ in
$\mathbb C^m$ such that $G$ is a trivial holomorphic vector bundle on $(\mathbb P^1-\{0,\infty\})\times X$
equipped with a meromorphic connection $\nabla$ with a regular singularity along $\infty\times X$, the Brieskorn lattice
$G_0$ is a trivial holomorphic vector bundle on
 $(\mathbb P^1-\{\infty\})\times X$ such that $G_0|_{(\mathbb P^1-\{0,\infty\})\times X}=G$, and the connection $\nabla$ has Poincar\'e
rank $\leq 1$ on $G_0$ along $0\times X$. When restricted to the
parameter $x=0$, this Brieskorn lattice coincides with the one
defined algebraically for the trivial deformation of $f$. Since we
can solve the Birkhoff problem for the trivial deformation, the
solution can be extended to a solution of the Birkhoff problem for
the Brieskorn lattice of the universal unfolding of $f$. We thus get
a Frobenius type structure on  the universal unfolding.

To get a Frobenius manifold
structure on the universal unfolding parameter space, one need to find a primitive form to transplant the
Frobenius type structure to the tangent sheaf of the universal
unfolding parameter space.  Another approach is to start with the solution of the Birkhoff
problem for a subdiagram deformation satisfying certain conditions, and then use a theorem of
Hertling and Manin \cite{HM} to show that this solution has a
universal deformation, which gives the Frobenius manifold structure on the  universal unfolding parameter space. This is
the Landau-Ginzburg B-model for the Laurent polynomial $f$.

\medskip
In summary, we start with the Brieskorn lattice for $f$, which is
obtained as the Fourier transform of the Gauss-Manin system
associated to $f:{\mathbb G}_m^n\to \mathbb A^1$.  Solve the
Birkhoff problem for it using Saito's criterion. The
Brieskorn lattice for $f$ have a deformation, which is the Brieskorn
lattice for the universal unfolding.  Extend the solution of
Birkhoff problem for the Brieskorn lattice associated to $f$ to the
solution of the Birkhoff problem for the Brieskorn lattice of the
universal unfolding. Or we can start with the Brieskorn lattice for a subdiagram
deformation. Solve the Birkhoff problem. Then apply the extension theorem of Hertling and Manin.

\medskip
Arithmetically, over a finite field $\mathbb F_q$ with $q$ elements
of characteristic $p$, we work with the Deligne-Fourier transform
$\mathcal F(Rf_!\overline{\mathbb Q}_\ell)$ of
$Rf_!\overline{\mathbb Q}_\ell$, where $\ell$ is a prime number
distinct from $p$. It should satisfy arithmetic counterpart of the
conditions for applying Saito's criterion, that is, it a lisse pure
sheaf on ${\mathbb P}^1-\{0,\infty\}$.  It is tamely ramified at
$0$, and its slopes at $\infty$ are $\leq 1$. Actually we can prove
such kind of results for $\mathcal F(RF_ {!}\overline {\mathbb
Q}_\ell)$ for any non-degenerate deformation $F$ of $f$ which
preserves the Newton polytope at $\infty$.

\section{$\ell$-adic realization of the Frobenius type structures}

In this section, we work over a finite ground field $\mathbb F_q$ with $q$ elements of characteristic $p$. Let
$\ell$ be a prime number distinct from $p$. Let $\psi:\mathbb F_q\to\overline{\mathbb Q}_\ell^\ast$ a
nontrivial additive character. For any $\mathbb F_q$-scheme $S$ of finite type, let $D_c^b(S, \overline{\mathbb Q}_\ell)$ be
the derived category of $\overline{\mathbb Q}_\ell$-schemes on $S$ defined in \cite[1.1.2]{D}. For any vector bundle
$E\to S$ of rank $r$, let $E^\vee$ be the dual vector bundle. The \emph{Deligne-Fourier transform} is the functor
$$\mathcal F_\psi: D_c^b(E,\overline{\mathbb Q}_\ell)\to D_c^b(E^\vee,\overline{\mathbb Q}_\ell),\quad
K\mapsto R\mathrm{pr}^\vee_!(\mathrm{pr}^\ast K\otimes
\langle\,,\,\rangle^\ast \mathcal L_\psi)[r],$$ where $E^\vee\to S$ is the dual vector bundle of $E$,  $\mathrm{pr}:E\times_S E^\vee\to E$
and $\mathrm{pr}^\vee:E\times_SE^\vee\to E^\vee$ are the projections, and $\langle\,,\,\rangle:E\times_SE^\vee\to\mathbb A_S^1$ is the pairing,
and $\mathcal L_\psi$ is the Artin-Schreier sheaf associated to the nontrivial additive character $\psi$.

Let $\mathbf{w}_j=(w_{1j},\ldots, w_{nj})\in\mathbb Z^n$
$(j=1,\ldots, N)$, let $f=\sum_{j=1}^N a_j t_1^{w_{1j}}\cdots
t_n^{w_{nj}}$  be a Laurent polynomial with nonzero coefficients
$a_j\in \mathbb F_q$,  and let $\Delta$ be the Newton polyhedron of
$f$ at $\infty$. Assume $0$ lies in the interior of $\Delta$, and
assume $f$ is non-degenerate, that is, for any face $\sigma$ of
$\Delta$ not containing $0$, the equations
$$\frac{\partial f_\sigma}{\partial t_1}=\cdots=\frac{\partial f_\sigma}{\partial
t_n}=0$$ define an empty subscheme in $\mathbb
G_{m,\mathbb F_q}^n=\mathrm{Spec}\,{{\mathbb F}_q}[t_1^{\pm 1},\ldots, t_n^{\pm 1}]$,
where $f_\sigma=\sum_{\mathbf w_j\in \sigma} a_jt_1^{w_{1j}}\cdots
t_n^{w_{nj}}$.

Let $g_1(t_1,\ldots, t_n), \ldots, g_m(t_1,\ldots, t_n)$ be a
family of Laurent polynomials. Consider the deformation
\begin{eqnarray*}
F_x(t)=f(t_1,\ldots, t_n)+x_1 g_1(t_1,\ldots, t_n)+\cdots + x_m
g_m(t_1,\ldots, t_n)
\end{eqnarray*} of $f$.
Suppose the Newton polytopes of $g_k$ $(k=1,\ldots, m)$ are
contained in $\Delta$. There exists a Zariski open subset $X\subset
\mathbb A^m$ containing the origin so that for any
$\overline{\mathbb F}_q$-point $(x_1,\ldots,x_m)$ in $X$, the Newton
polyhedron of $F_x(t)$ at $\infty$ coincides with $\Delta$, and
$F_x(t)$ is non-degenerate with respect to $\Delta$. Consider the
morphism
$$F: {\mathbb G}_{m,X}^n\to \mathbb A_X^1,\quad (t,x)\mapsto
(F_x(t),x).$$ The Deligne-Fourier transform $\mathcal
F_\psi(RF_!\overline{\mathbb Q}_\ell)$ of $RF_!\overline{\mathbb
Q}_\ell$ is an analogue of $G$ in \S 2. Here the Deligne-Fourier
transform is taken for the vector bundle $\mathbb A_X^1\to X$. Note
that we have
$$\mathcal F_{\psi^{-1}}(RF_!\overline{\mathbb Q}_\ell)\cong a^\ast
\mathcal F_\psi(RF_!\overline{\mathbb Q}_\ell),$$ where $a:\mathbb
A_X^1\to\mathbb A_X^1$ is the morphism $\tau\mapsto -\tau$.

\begin{theorem}\label{mainthm}  Notation as above. Suppose the Newton polytopes of
$g_i$ are contained in $\Delta$.

(i)  When restricted to  $(\mathbb P^1-\{0,\infty\})\times X$,
$\mathscr H^i\Big(\mathcal F_\psi(RF_!\overline{\mathbb
Q}_\ell)\Big)=0$ for $i\not=n-1$, and $\mathcal H^{n-1}\Big(\mathcal
F_\psi(RF_!\overline{\mathbb Q}_\ell)\Big)$ is a pure lisse sheaf of
weight $n$.

(ii) When restricted to  $\mathbb P_X^1-(\{0,\infty\}\times X)$, we
have a perfect pairing
$$\mathcal H^{n-1}\Big(\mathcal
F_\psi(RF_!\overline{\mathbb Q}_\ell)\Big) \times \mathcal
H^{n-1}\Big(\mathcal F_{\psi^{-1}}(RF_!\overline{\mathbb
Q}_\ell)\Big) \to \overline{\mathbb Q}_\ell(-n).$$

(iii) For each fixed geometric point $x$ of $X$, the restriction of
$\mathcal H^{n-1}\Big(\mathcal F_\psi(RF_!\overline{\mathbb
Q}_\ell)\Big)$ to $(\mathbb P^1-\{0,\infty\})\times \{x\}$ is tamely
ramified at $0$, and has slopes $\leq 1$ at $\infty$.
\end{theorem}

Let $K$ be a local field with perfect residue field. We have an
equivalence between the category of $\overline{\mathbb
Q}_\ell$-sheaves on $\mathrm{Spec}\, K$ and the category of
$\overline{\mathbb Q}_\ell$-representations of
$\mathrm{Gal}(\overline K/K)$. The higher ramification subgroups of
$\mathrm{Gal}(\overline K/K)$ give rise to a decreasing filtration
in upper numbering $\{G^{(r)}\}_{r\in {\mathbb Q}_{\geq 0}}$. Confer
\cite{S}. Let $G^{(r+)}$ be the closure of $\cup_{s>r}G^{(s)}$. Then
$G^{(0+)}$ is the wild ramification subgroup of
$\mathrm{Gal}(\overline K/K)$. Any $\overline{\mathbb
Q}_\ell$-representation $V$ of $\mathrm{Gal}(\overline K/K)$ is
semisimple as a representation of $G^{(0+)}$. Let $V=\oplus_i V_i$
be a decomposition of $V$ into irreducible representations of
$G^{(0+)}$. The slope of $V_i$ is the smallest rational number $s_i$
such that $V_i^{G^{(s_i+)}}=V$. The numbers $s_i$ are called slopes
of $V$.

Let $\mathbb F$ be the algebraic closure of $\mathbb F_q$, and let
$\mathbb F[[\pi]]$ be the formal power series ring. The field of
fractions of $\mathbb F[[\pi]]$ is the field $\mathbb F((\pi))$ of
formal Laurent series. On $\mathrm{Spec}\,\mathbb F[[\pi]]$, let $o$
be the divisor defined by the maximal ideal of of $\mathbb
F[[\pi]]$. Let $\mathcal H$ be a lisse sheaf on $(\mathbb
P^1-\{0,\infty\})\times X$. We say $\mathcal H$ is \emph{tamely
ramified} at $0\times X$ if for any $\mathbb F_q$-morphism
$g:\mathrm{Spec}\, \mathbb F[[\pi]]\to \mathbb P^1\times X$ such
that $g^\ast(0\times X)$ is the divisor $eo$ $(e\geq 1)$, the sheaf
$(g^\ast \mathcal H)|_{\mathrm{Spec}\, \mathbb F((\pi))}$ is tamely
ramified. We say $\mathcal H$ has \emph{slope $\leq r$} at
$\infty\times X$ if for any $\mathbb F_q$-morphism
$g:\mathrm{Spec}\, \mathbb F[[\pi]]\to \mathbb P^1\times X$ such
that $g^\ast(\infty\times X)$ is the divisor $eo$ $(e\geq 1)$, the
sheaf $(g^\ast \mathcal H)|_{\mathrm{Spec} \mathbb F((\pi))}$ has
slopes $\leq er$. We can also define slopes using Abbes-Saito's
theory for higher ramifications of Galois representations of local
field with imperfect residue field.

In view of the fact that a Frobenius type structure has logarithmic
pole at $\infty\times X$, and has Poincar\'e rank $\leq 1$ at
$0\times X$, the following fact which is more general than Theorem
\ref{mainthm}(iii) should be true: $\mathcal H^{n-1}\Big(\mathcal
F_\psi(RF_!\overline{\mathbb Q}_\ell)\Big)$ is tamely ramified at
$0\times X$, and has slopes $\leq 1$ at $\infty\times X$.

We have seen that connections of Poincar\'e rank $\leq 1$ gives rise
to structures such as Higgs fields. Due to our lack of explicit
description of the higher ramifications, we haven't been able to
extract structures hidden in $\ell$-adic sheaves of slope $\leq 1$
along a divisor.

We define an \emph{$\ell$-adic Frobenius type structure with a
metric} to be a pure lisse sheaf $\mathcal H$ on $(\mathbb
P^1-\{0,\infty\})\times X$ which is tamely ramified at $\infty\times
X$, and has slopes $\leq 1$ at $0\times X$, and has a perfect
pairing
$$\mathcal H\times a^\ast\mathcal H\to \overline{\mathbb
Q}_\ell(-n),$$ where $n$ is the weight of $\mathcal H$. By the above
discussion, $\mathrm{inv}^\ast\Big(\mathcal H^{n-1}\Big(\mathcal
F_\psi(RF_!\overline{\mathbb Q}_\ell)\Big)\Big)$ should define an
$\ell$-adic Frobenius type structure, where $\mathrm{inv}:(\mathbb
P^1-\{0,\infty\})\times X\to (\mathbb P^1-\{0,\infty\})\times X$ is
the morphism defined by $\tau\mapsto \frac{1}{\tau}$.

\medskip
The proof of Theorem \ref{mainthm} uses the properties of
$\ell$-adic Gelfand-Kapranov-Zelevinsky (GKZ) hypergeometric sheaves
introduced in \cite{F}. Choose $\mathbf w_{N+1},\ldots, \mathbf
w_{N'}$ so that
$$\Delta\cap \mathbb Z^n=\{\mathbf w_1,\ldots, \mathbf w_N,\mathbf
w_{N+1}, \ldots, \mathbf w_{N'}\}.$$ Let $\pi_2: \mathbb
G_m^n\times\mathbb A^{N'}\to \mathbb A^{N'}$ be the projection, and
let $H$ be the morphism $$H:\mathbb G_m^n\times \mathbb A^{N'}\to
\mathbb A^1,\quad (t_1,\ldots, t_n, y_1,\ldots, y_{N'})\mapsto
\sum_{j=1}^{N'}y_jt_1^{w_{1j}}\cdots t_n^{w_{nj}}.$$ We define the
\emph{$\ell$-adic GKZ hypergeometric sheaf}  to be the object in
$D_c^b(\mathbb A^{N'}, \overline{\mathbb Q}_\ell)$ defined by
$$\mathrm{Hyp}_\psi=R\pi_{2!}H^\ast \mathscr
L_\psi[n+N'].$$  We have the following:

\begin{theorem} \label{gkz} ${}$

(i) $\mathrm{Hyp}_\psi$ is a pure perverse sheaf on $\mathbb A^{N'}$
of weight $n+N'$ and of rank $(-1)^{N'}n!\mathrm{vol}(\Delta)$.

(ii) Suppose $V$ is a Zariski open subset of $\mathbb A^{N'}$ such
that for any $(a_1,\ldots, a_{N'})\in V(\overline {\mathbb F}_q)$,
the Laurent polynomial $\sum_{j=1}^{N'} a_j t_1^{w_{1j}}\cdots
t_n^{w_{nj}}$ is nondegenerate with respect to $\Delta$. Then
$\mathcal H^i(\mathrm{Hyp})|_V=0$ for $i\not=-N'$, and $\mathcal
H^{-N'}(\mathrm{Hyp}_\psi)|_V$ is lisse, pure of weight $n$, and of
rank $n!\mathrm{vol}(\Delta)$.

(iii) The Verdier dual $D(\mathrm{Hyp}_\psi):=R\mathcal
Hom(\mathrm{Hyp}_\psi,\overline{\mathbb Q}_\ell(N')[2N'])$ of
$\mathrm{Hyp}_\psi$ is isomorphic to
$\mathrm{Hyp}_{\psi^{-1}}(n+N')$, where $(n+N')$ denotes the Tate
twist by $\overline{\mathbb Q}_\ell(n+N')$. In particular, on the
open set $V$, we have a perfect pairing
$$\mathcal H^{-N'}(\mathrm{Hyp}_\psi)\times \mathcal
H^{-N'}(\mathrm{Hyp}_{\psi^{-1}})\to \overline{\mathbb
Q}_\ell(-n).$$
\end{theorem}

The main technique to study the $\ell$-adic GKZ hypergeometric sheaf
is again the Deligne-Fourier transform. Let $\iota:\mathbb T^n\to
\mathbb A^{N'}$ be the morphism defined by
$$(t_1,\ldots, t_n)\mapsto (t_1^{w_{1j}}\cdots t_n^{w_{nj}})_{j=1,\ldots,
N'}.$$ In \cite{F}, we prove that $$\mathrm{Hyp}_\psi=\mathcal
F_\psi(\iota_!\overline{\mathbb Q}_\ell[n]),$$ where $\mathcal
F_\psi$ is the Deligne-Fourier transform for the vector bundle
$\mathbb A^{N'}\to\mathrm{Spec}\,\mathbb F_q$. Using the assumption
that $0$ lies in the interior $\Delta$, one can show
$\iota_!\overline{\mathbb Q}_\ell[n]=\iota_{!\ast}\overline{\mathbb
Q}_\ell[n]$. From the standard facts on perverse sheaves and the
Deligne-Fourier transform, one deduces that $\mathrm{Hyp}_\psi$ is a
pure perverse sheaf of weight $n+N'$, and $D(\mathrm{Hyp}_\psi)\cong
\mathrm{Hyp}_{\psi^{-1}}(n+N')$. The other statements in Theorem
\ref{gkz} require a detailed study of the morphism $H$ relatively to
the toric compactification defined by the convex polytope $\Delta$.

\medskip
Write $$F_x(t)=f(t_1,\ldots, t_n)+x_1 g_1(t_1,\ldots, t_n)+\cdots + x_m
g_m(t_1,\ldots, t_n)=\sum_{j=1}^{N'} y_j(x_1,\ldots, x_m) t_1^{w_{1j}}\cdots t_n^{w_{nj}}.$$
Then $y_j(x_1,\ldots,x_m)$ are linear polynomial of $x_1,\ldots,x_m$ with coefficients depending on
$f$ and $g_1,\ldots, g_m$.  Consider the morphism
$$\Phi:\mathbb A_X^1\to \mathbb A^{N'}, \quad (\tau, x_1,\ldots, x_m)\to
(\tau y_1(x_1,\ldots, x_m), \ldots,\tau y_{N'}(x_1,\ldots, x_m) ).$$
Theorem \ref{mainthm} follows from Theorem \ref{gkz} and the
following:

\begin{proposition} \label{link} ${}$

(i) We have $\mathcal F_\psi(RF_!\overline{\mathbb Q}_\ell)\cong
\Phi^\ast (\mathrm{Hyp}_\psi)[1-(n+N')]$.

(ii) The image of $(\mathbb P^1-\{0,\infty\})\times X$ under $\Phi$
is contained in the set $V$ parameterizing non-degenerate Laurent
polynomials.
\end{proposition}

Instead of working with Laurent polynomials, one can also work with polynomials, and similar results still hold.

\medskip
Originally our motivation for studying the $\ell$-adic GKZ
hypergeometric sheaf comes from the study of exponential sums using
$\ell$-adic cohomology theory. For any $\mathbb F_q$-rations points
$x=(x_1,\ldots, x_{N'})$ of $\mathbb A^{N'}$, it follows from the
Grothendieck trace formula that we have
\begin{eqnarray}
(-1)^{n+N}\mathrm{Tr}(\mathrm{Fr}_x, (\mathrm{Hyp}_\psi)_{\bar
x})=\sum_{t_1,\ldots, t_n \in\mathbb F_q^\ast}\psi(\sum_{j=1}^{N'}
y_j t_1^{w_{1j}}\cdots t_n^{w_{nj}}),
\end{eqnarray} where the left hand side is the trace of the
geometric Frobenius at the point $x$ acting on the stalk of
$\mathrm{Hyp}_\psi$ at the geometric point $\bar x$ above $x$. Note
that the right hand side of the above equation is a family of
exponential sums parameterized by $(y_1,\ldots, y_{N'})$. It is an
analogue of the oscillatory integral
$$\int_{\sigma} e^{i\sum_{j=1}^{N'} y_j t_1^{w_{1j}}\cdots t_n^{w_{nj}}}
\frac{dt_1}{t_1}\cdots \frac{dt_n}{t_n},$$ where $\sigma$ is an
$n$-dimensional cycle in $\mathbb G_{m,\mathbb C}^n$. This integral
is a solution of the GKZ hypergeometric system of differential
equations.

Similarly, for any $\mathbb F_q$rational point $(\tau,x)=(\tau,
x_1,\ldots, x_m)$ of $\mathbb A^1_X$, we have
\begin{eqnarray}
\begin{array}{l}
\quad -\mathrm{Tr}\Big(\mathrm{Fr}_{(\tau,x)}, (\mathcal
F_\psi(RF_!\overline {\mathbb
Q}_\ell))_{\overline{(\tau,x)}}\Big)\\
=\sum_{t_1,\ldots, t_n\in \mathbb F_q^\ast}\psi\Big(\tau
(f(t_1,\ldots, t_n)+x_1 g_1(t_1,\ldots, t_n)+\cdots + x_m
g_m(t_1,\ldots, t_n))\Big).
\end{array}
\end{eqnarray}
It is clear that the family of exponential sum on the right hand
side of (3) is the composite of $\Phi$ and the family of exponential
sum on the right hand side of (2). This gives an explanation of
Proposition \ref{link} (i) on the level of functions. The
oscillatory integral corresponding to the exponential sum in (3) is
$$\int_\sigma e^{i\tau
(f(t_1,\ldots, t_n)+x_1 g_1(t_1,\ldots, t_n)+\cdots + x_m
g_m(t_1,\ldots, t_n))}\frac{dt_1}{t_1}\cdots \frac{dt_n}{t_n}.$$ It
is a solution of the system of differential equations of the
$D$-module defined by the connection introduced in \S 2.

\end{document}